\documentclass[preprint,11pt]{article}
\usepackage{amssymb}
\usepackage{amsmath}
\usepackage{amsthm}
\usepackage{enumerate}
\usepackage{graphicx}
\usepackage{amsfonts}
\usepackage{authblk}
\usepackage{mathtools,xcolor,url, todonotes}

\newtheorem{example}{Example}

\newtheorem{problem}{Problem}[section]

\newcommand {\RR} {\mathbb R}

\makeatletter
\def\dddots{\mathinner{\mkern1mu\raise\p@
    \hbox{.}\mkern2mu\raise4\p@\hbox{.}\mkern2mu
    \raise7\p@\vbox{\kern7\p@\hbox{.}}\mkern1mu}}%
\makeatother
\usepackage[english]{babel}

\title{Polynomial mapped bases: theory and applications}
	
\author[a,e]{S. De Marchi}
\author[b,e]{G. Elefante}
\author[c,e]{E. Francomano}
\author[d,e]{F. Marchetti}
\affil[a]{University of Padova, Department of Mathematics “Tullio Levi-Civita”, Italy,  \texttt{stefano.demarchi@unipd.it}}
\affil[b]{University of Padova, Department of Mathematics “Tullio Levi-Civita”, Italy,  \texttt{giacomo.elefante@unipd.it}}
\affil[c]{University of Palermo, Department of Engineering, Italy, \texttt{elisa.francomano@unipa.it}}
\affil[d]{University of Padova, Department of Mathematics “Tullio Levi-Civita”, Italy,  \texttt{francesco.marchetti@unipd.it}}
\affil[e]{IN$\delta$AM Gruppo Nazionale di Calcolo Scientifico}

\begin{document}

\maketitle

\begin{abstract} 	
In this paper, we collect the basic theory and the most important applications of a novel technique that has shown to be suitable for scattered data interpolation, quadrature, bio-imaging reconstruction. The method relies on \textit{polynomial mapped bases} allowing, for instance, to incorporate data or function discontinuities in a suitable \emph{mapping} function. The new technique substantially mitigates the Runge's and Gibbs effects.
\end{abstract}

\section{Introduction}
The reconstruction of functions from data is a historical and common task in many applications, as well detailed in the fundamental books \cite{Davis, deBoor,T12}.

It is well known that the {\it Runge's effect} may arise in the approximation of particular functions with certain nodes distributions \cite{Runge}. Moreover, when discontinuities or jumps occur, the reconstruction task is even more interesting and challenging due to the well-known {\it Gibbs phenomenon}, which appears especially near the discontinuities \cite{gibbs}. 

For example, both phenomena are of interest in the context of medical image analysis (such as in Computerized Tomography (CT), in Magnetic Resonance (MR), and their variants (SPECT, fMRI), in Magnetic Particle Imaging (MPI)), where the images often need to be geometrically aligned, registered or simply reconstructed by a proper (re)sampling. In bio-imaging the Gibbs effects are often referred to as \textit{ringing artifacts}. General approaches to bypass unavoidable reconstruction instabilities are as follows:
\begin{enumerate}
\item by a clever choice of interpolation points (cf. e.g. \cite{DeEM21}), 
\item by using rational approximation (cf. e.g. \cite{BDEM21,SarraBai18}), 
\item by filtering techniques (cf. e.g. \cite{DEM17}). 
\end{enumerate}
In the last years, we developed a new simple technique that we called {\it fake nodes} approach, which is a {\it polynomial mapped bases} method \cite{4,5}. This paper collects the most important results that the CAA Research Group \cite{2} and their collaborators achieved about this approach. In Section 2, we outline the theoretical background behind the mapped bases approach. In Section 3 we describe the application to numerical quadrature, and in Section 4 we provide a list of some open problems on which we are still working. Finally, we conclude in Section 5.

\section{Mapped bases approach} \label{sec:2}

Let $X_{N}=\{\boldsymbol{x}_i, \; i=1, \ldots,N \} \subset \Omega$, $\Omega \subseteq \mathbb{R}^d$ be a set of distinct nodes and let $F_{N}=\{f(\boldsymbol{x}_i), \; i=1, \ldots,N \}$ be the function values at $X_N$, which are sampled from a function $f: \Omega \longrightarrow \mathbb{R}$. Letting ${\cal B}_N$ be a given finite dimensional function space, the multivariate interpolation problem consists in finding a function $P_f \in {\cal B}_N$ such that
\begin{equation*}
    P_f(\boldsymbol{x}_i) = f(\boldsymbol{x}_i), \quad i=1, \ldots, N. 
\end{equation*}
We assume $P_f \in {\cal B}_N \coloneqq \textrm{span}\{B_1, \ldots, B_N\},$ where $B_i: \Omega \longrightarrow \mathbb{R}$, $i=1,\ldots, N$, are the basis functions. 

In approximation practice, when samples are given, {\it resampling} is often necessary. This can be done by choosing {\it good} interpolation points (for instance Chebyshev points). Moreover, depending on applications we can extract {\it mock-Chebyshev} points from equispaced samples, Approximate Fekete points, Discrete Leja Sequences \cite{BDSV10} or (P, $f$, $\beta$)-greedy points \cite{WSH21} when the data are scattered. 
When the underlying function presents a steep gradient in one or more locations, in order to achieve a faster convergence, the nodes could be mapped via a conformal map which allows to cluster them in a precise location, as in \cite{BerEle}. Recently in \cite{AP16}, the authors investigated a weighted least-squares approximation via mapped polynomial bases of the interval $[-1,1]$, by using the so-called \textit{Kosloff Tal-Ezer map} \cite{Kosloff},
$$T_\alpha(x)=\frac{ \sin(c_\alpha x)}{\sin(c_\alpha)},\; x \in [-1,1], \; \alpha \in (0,1], \; c_\alpha=\alpha \, \frac{\pi}{ 2}\,.$$
This gives rise to the $\alpha$-polynomial space
$\mathbb{P}_n^\alpha=\{ p \circ T_\alpha, \; p \in \mathbb{P}_n \}\,$  
which corresponds to the space of trigonometric polynomials when $\alpha=1$.

In our mapped bases approach, we consider an injective map $S: \Omega \subset \mathbb{R}^d \longrightarrow \mathbb{R}^d$ and the mapped basis $\{B_i^S,\;i=1,2,\ldots,N\}$ with $B_i^S=B_i\circ S$. Then, we construct the interpolant $R_f \in {\cal B}^S_N \coloneqq \textrm{span}\{B^S_1, \ldots, B^S_N\}$ of the function $f$ as
\begin{equation*}
R_f(\boldsymbol{x})=  \sum_{i=1}^N \alpha^S_i B^S_i(\boldsymbol{x}) = \sum_{i=1}^N \alpha^S_i B_i(S(\boldsymbol{x}))=P_g(S(\boldsymbol{x})), \;\; \forall \, \boldsymbol{x} \in \Omega\,.
\end{equation*}
The function $g$ has the {\it no-resampling property} that is $g_{|S(X_N)}=f_{|X_N}$ (cf. \cite{4}). Thus, having the mapped basis ${\cal B}^S_N$, the construction of the interpolant $R_f$ is then equivalent to build the classical interpolant $P_g \in {\cal B}_N$ at the mapped nodes $S(X_N)$, which we called \textit{fake nodes}. The term {\it fake} has been introduced, because the new points are employed only in the final reconstruction step.

Here we summarize the most important properties of the aforementioned approach, while we refer to \cite{4,5} for a detailed presentation.
\begin{itemize}
\item {\it Generality}. The mapped bases approach can be applied to any basis spanning the approximation space.
\item {\it Equivalence of the Lebesgue constants}. The Lebesgue constant \cite{Brutman} of the points mapped via $R_f^S$ is equivalent to that of the image $\Omega$ through $S$.  Indeed, the Lebesgue constant $\Lambda^S(\Omega)$ associated to the mapped nodes satisfies the equation $$\Lambda^S(\Omega) = \Lambda(S(\Omega)).$$
\item {\it Stability of the mapped interpolant}.  Let $\boldsymbol{f}$ be the associated vector of function values and $\boldsymbol{\tilde f}$ be the vector of perturbed values. Let $R^S_{{f}}$ and $R^S_{\tilde{f}}$ be the interpolants of the function values $\boldsymbol{f}$ and $\boldsymbol{\tilde{f}}$ respectively. Then,
\begin{equation*}
    \lvert\lvert R_f^S - R_{\tilde f}^S \rvert\rvert_{\infty,\Omega} \leq  \Lambda^S(\Omega) \;  \| \boldsymbol{f}-\boldsymbol{\tilde{f}} \|_{\infty, X_N}.
\end{equation*}
\item {\it Error bound inheritance}. For any given
function norm, we have
$$\lvert\lvert R_f^S-f  \rvert\rvert_{\Omega}  = \lvert\lvert P_g-g   \rvert\rvert_{S(\Omega)}, $$
where $ f = g\circ S$. 

\end{itemize}
\begin{problem}\label{prob2}
How can we find a suitable map $S$ for mitigating the Runge's and Gibbs effects?  
\end{problem}

\subsection{Polynomial mapped bases}
We now describe the polynomial mapped bases approach firstly in the univariate case and secondly in the multivariate case.
\subsubsection{The univariate case}

Let $K=[a,b] \subset \RR$ and let ${\cal M}=\{1,x,x^2,\ldots,x^n\}$ be the basis of monomials of the space of the univariate polynomials $\mathbb{P}_n(\RR)$. If $N=n+1$, the univariate polynomial interpolation problem at $X_N$ has unique solution if and only if the Vandermonde determinant
\begin{equation*}\label{Vdm1}
\text{VDM}(X_N;{\cal M})=\prod_{i <j} (x_i-x_j)
\end{equation*} 
is so that $\text{VDM}(X_N;{\cal M})\ne 0$. We point out that this implies $\text{VDM}(X_N;S({\cal M}))\ne 0$, i.e., the well-posedness of the mapped bases scheme.

Two simple algorithms, {$S$-Runge} and $S$-Gibbs, have been designed to provide a constructive solution to Problem \ref{prob2} (for details see \cite[Algorithms 1 \& 2]{4}). In the $S$-Runge, the map $S$ is basically a linear transformation from the set $X_N$ of equispaced or random points to the Chebyshev-Lobatto (CL) ones, which represents a ``natural'' choice for a stable interpolation. For {$S$-Gibbs}, first we need to identify the set 
$${D}_{m}:=\left\{(\xi_i,d_i)\:|\:\xi_i\in K\setminus \partial K,\;\xi_i<\xi_{i+1},\;\text{ and }d_i\coloneqq |f(\xi_i^+)-f(\xi_i^-)|,\; i=1,\dots,m\right\}$$
of all $m$ discontinuities of $f$. Then, the map $S$-Gibbs is constructed to be a linear map depending on the set of discontinuities $D_m$ and on a {\it shift parameter} $k$, and its role consists in increasing the distance between the nodes just before and after the discontinuity locations. We experimentally observed that the tuning of the shift parameter is not critical, since the resulting interpolation process {is not sensitive to its choice, provided that it is sufficiently large}, i.e. in such a way that in the mapped space the so-constructed function $g$ has no steep gradients. 

A stable S-Gibbs method, called "Gibbs-Runge-Avoiding Stable Polynomial 
Approximation", shortly {\it GRASPA}, has been recently proposed in \cite{DeEM21}, where the Lebesgue constant has been related to the shift parameter $k$. In GRASPA we consider the partition $\mathcal{D}=\{K^1,\dots,K^{m+1}\}$ of $K$, in which each subinterval $K^i\subset K$ is separated to the following one by an element of $D_m$. Precisely, each $K^i$ contains a discontinuity of the underlying function. In Example \ref{esempiooo}, we show the effectiveness of the GRASPA scheme by considering a test function with one discontinuity.

In practical applications, we remark that the exact position of the discontinuities is not at disposal, and can be investigated by using well-known and stable techniques, such as the Canny edge-detection algorithm described in \cite{Canny} or, for irregularly samples signals and images, as presented in \cite{Archibald05}. When radial basis functions are used, the analysis of the coefficients of the interpolant can give information on the location of the discontinuities, as described in \cite{Romani19}. Recently we proposed another approach to extract the location of the discontinuities through a segmentation method based on a classification algorithm from machine learning  \cite{DEMPR20}.

\begin{example}\label{esempiooo}
Let us consider the function $f_1$ on $K_1=[-1,1]$
\begin{equation*}
    f_1(x)=
    \begin{dcases}
        -\bigg(\frac{1}{10((x+0.5)^2+0.1)}\bigg)^5 & \textrm{if $-1\le x\le 0$,}\\
        e^{-x} & \textrm{if $ 0<x\le 1$,}
    \end{dcases}
\end{equation*}
which is discontinuous at $\xi=0$. Therefore, $\mathcal{D}=\{K^1,K^2\}$ with $K^1=[-1,\xi],K^2=]\xi,1]$. In Figg. \ref{fig3} and \ref{fig3bis}, we compare the results achieved by different interpolation methods.
\begin{figure}[!hbt]
  \centering
  \includegraphics[width=0.32\linewidth]{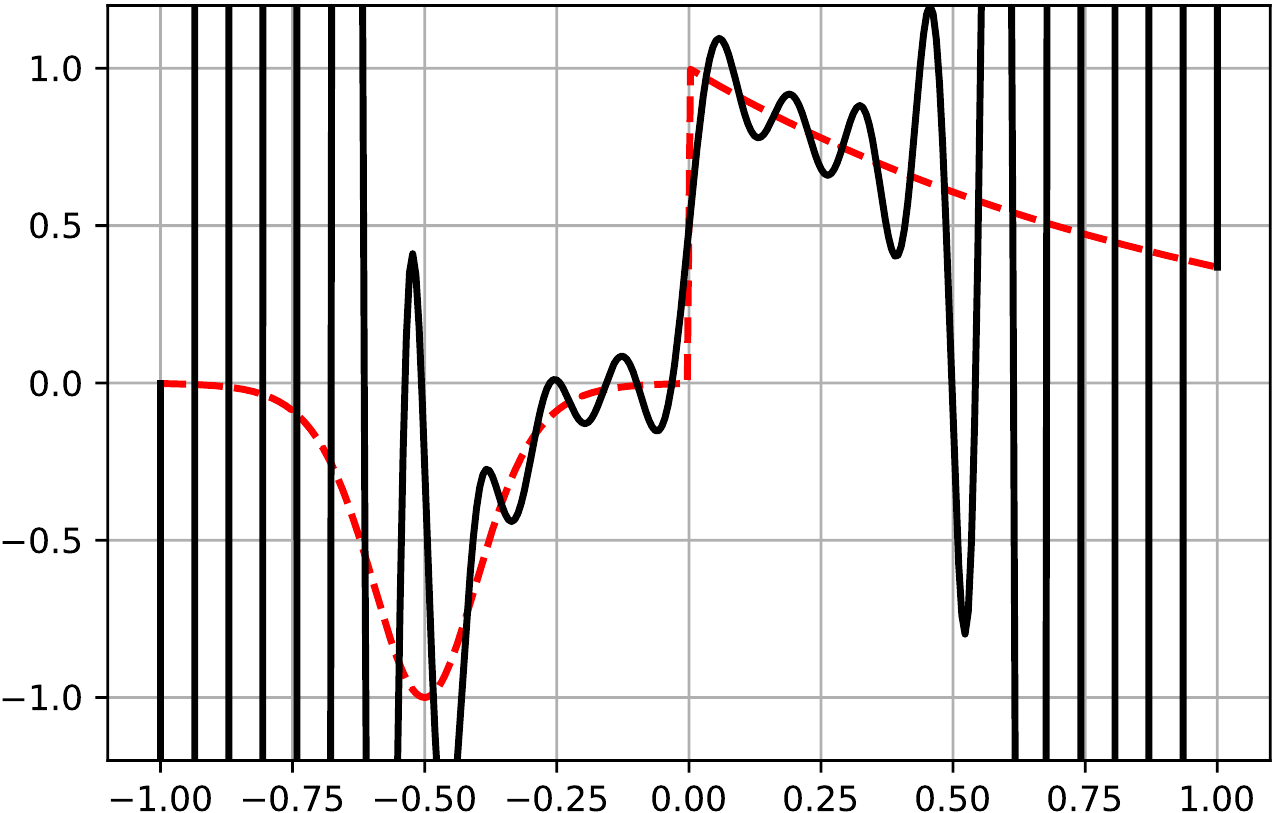}  
  \includegraphics[width=0.32\linewidth]{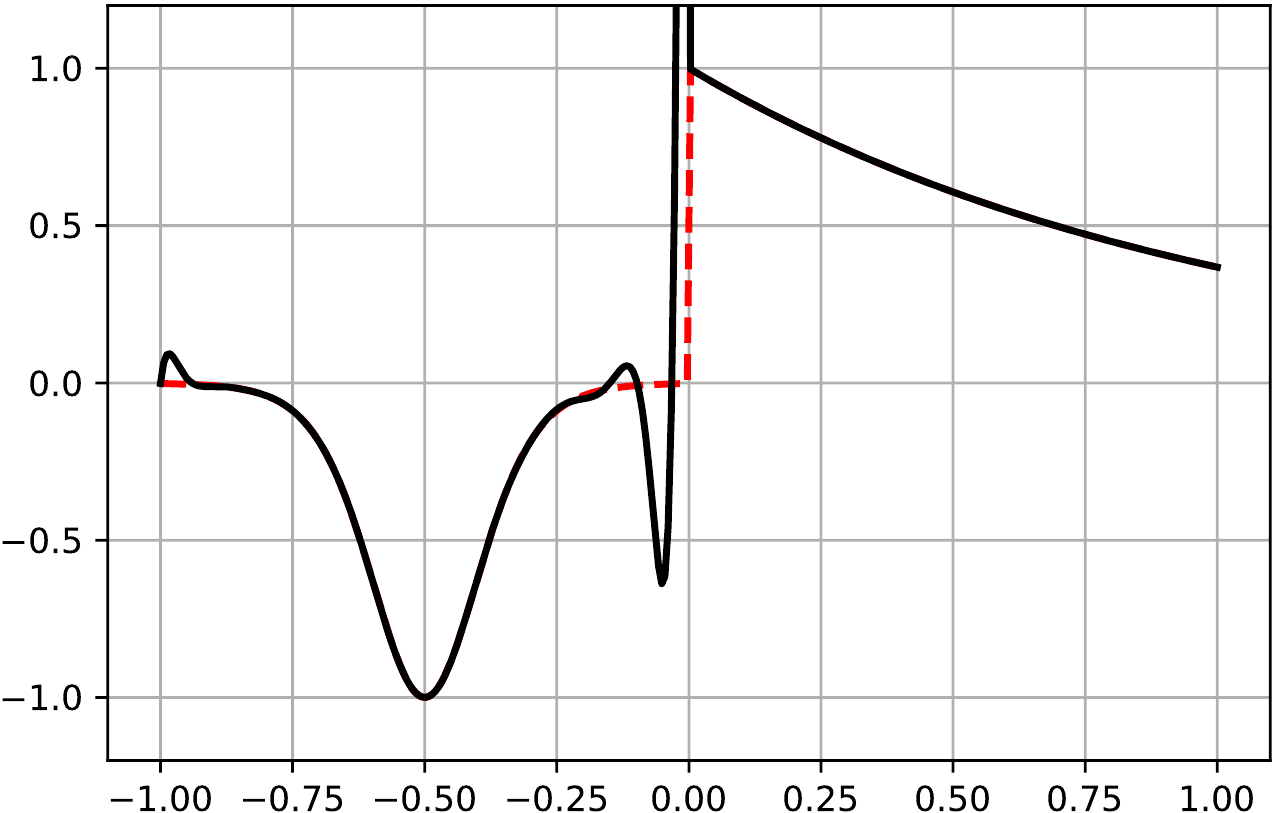}
   \includegraphics[width=0.32\linewidth]{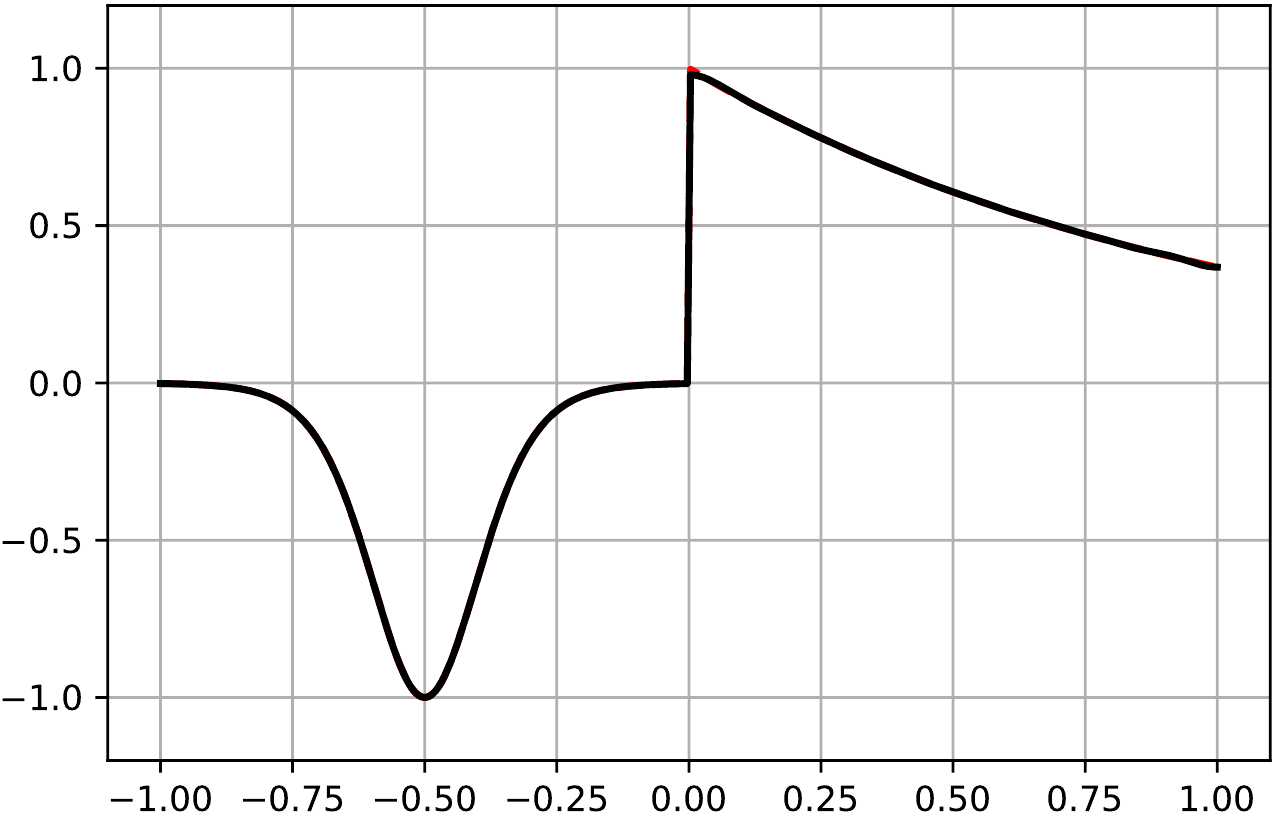}  
\caption{The function $f_1$ in dashed red and the interpolant with $n=32$ equispaced nodes in black. From left to right: classical, S-Gibbs and GRASPA approach, respectively.}
\label{fig3}
\end{figure}

\begin{figure}[!hbt]
  \centering
  \includegraphics[width=0.45\linewidth]{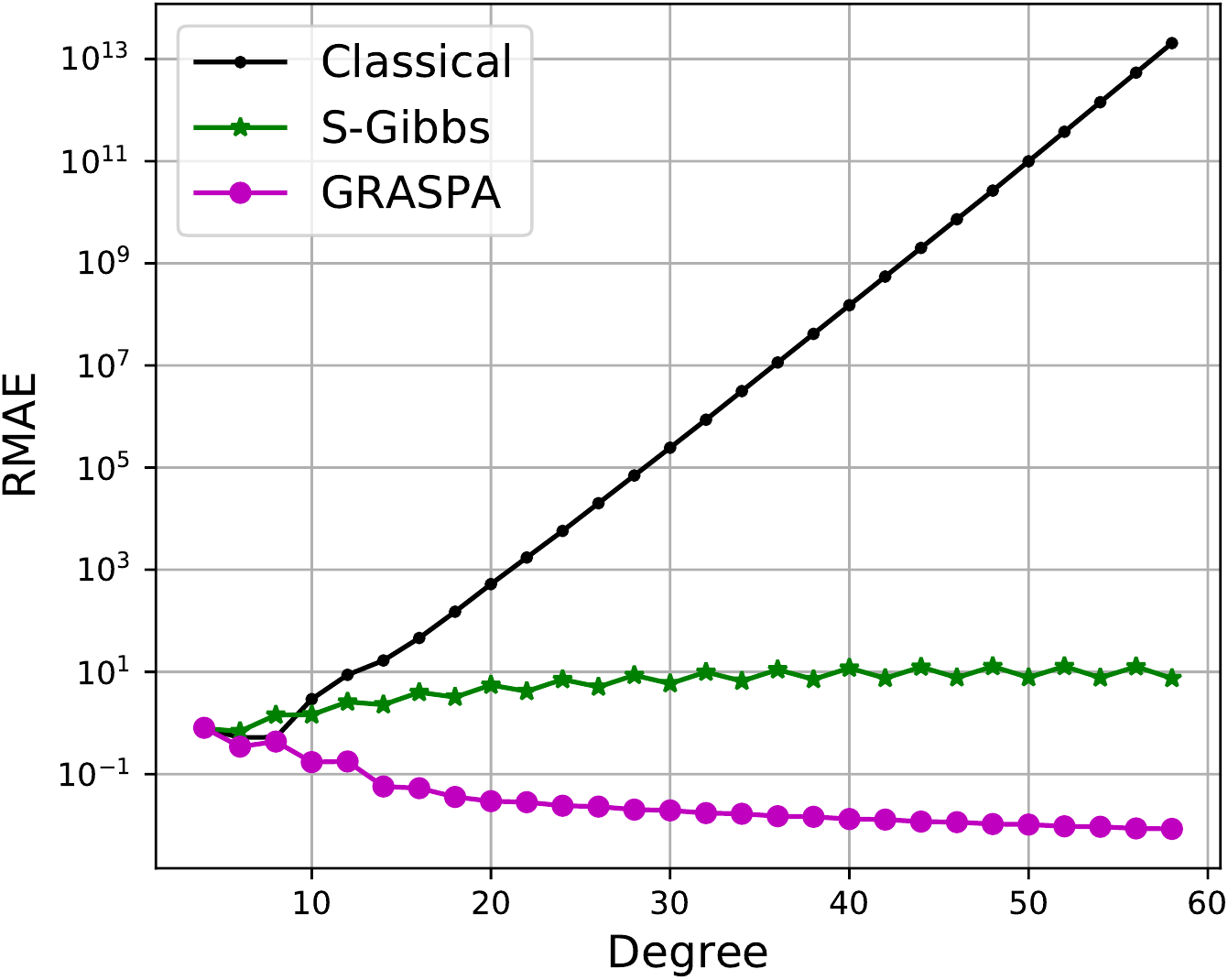}  
\caption{Relative Maximum Absolute Errors (RMAE) related to the three approaches.}
\label{fig3bis}
\end{figure}

As we can notice, the S-Gibbs map resolves the Gibbs phenomenon. However, if the Runge's phenomenon takes place in some subdomain the interpolating function diverges, but by means of the GRASPA map we could prevent the appearance of both.
\end{example}

\subsubsection{The multidimensional case}

The definition of a S-Runge-like algorithm in any dimension $d>1$ still represents an open problem. Nevertheless, some results have already been proved. 

In the case $d=2$ in the square $[-1,1]^2$, as far as the Runge's phenomenon is concerned, the optimal points are explicitly known to be the {\it Padua points} ${\rm Pad}_n$ \cite{BDVX07,CDV05}. In \cite{5}, we proposed the use of the map
$ S: [-1,1]^2 \longrightarrow [-1,1]^2$ 
$$ S(\boldsymbol{x}) = \left(-\cos\bigg(\pi\frac{\boldsymbol{e}_1^\intercal\boldsymbol{x}+1}{2}\bigg), -\cos\bigg(\pi\frac{\boldsymbol{e}_2^\intercal\boldsymbol{x}+1}{2}\bigg)\right).$$
where $\boldsymbol{x}=(x_1,x_2)$ and $\boldsymbol{e}_i$, $i=1,2$, are the unit vectors of $\mathbb{R}^2.$ This function maps the set of nodes
\begin{equation}
\label{NodesRunge}
    {X}_N =  \left\lbrace \left(\frac{2(i-1)}{n}-1, \frac{2(j-1)}{n+1}-1\right), \begin{array}{c} i = 1,\dots,n+1\\j= 1, \dots, n+2 \end{array},  
    \begin{array}{c} i+j \equiv 0\\ \pmod{2} \end{array}   \right\rbrace,
\end{equation}
where $N=(n+1)(n+2)/2$, onto the Padua points. In higher dimensions, where Padua points are not known, we may proceed analogously by considering the so-called {Lissajous points} \cite{Erb15}.

Concerning the mitigation of the Gibbs phenomenon, a straightforward extension of the S-Gibbs algorithm to general dimensions $d>1$ has been proposed in \cite{5} and applied in \cite{DEFMPP20} to kernel-based approximation with discontinuous kernels in the framework of the MPI and in \cite{PCD22} in multimodal medical imaging. In this context, it is a common practice to undersample the anatomically-derived segmentation images to measure the mean activity of a co-acquired functional image. This avoids the resampling-related Gibbs effect that would occur in oversampling the functional image. It turns out that the mapped bases scheme provides a reduction of the Gibbs effect when oversampling the functional image, as proved by a tight error analysis (we refer to \cite{PCD22} for further details).

\subsection{Barycentric rational mapped bases}
As well-known, it is possible to write the unique polynomial $P_f$ of degree at most $n$ interpolating $f$ at the set $X_N\subset K$ in the second barycentric form
\begin{equation}
    P_f(x)= \frac{ \sum_{i=1}^N \frac{\lambda_i }{ x-x_i} f_i}{ \sum_{i=1}^N \frac{ \lambda_i }{ x-x_i}},\; x\in K,
\end{equation}
where $\lambda_i=\displaystyle \prod_{j\neq i} \frac{1}{x_i-x_j}$ are the so-called {\it weights} (see \cite{BerTre}).
If the $\lambda_i$, $i=1,\dots,N$, are changed to other nonzero weights, say $w_i$, then the corresponding barycentric rational interpolant is
\begin{equation}\label{BarInt}
    r_f(x)= \frac{\sum_{i=1}^N \frac{w_i}{x-x_i}f_i}{\sum_{i=1}^N \frac{w_i}{x-x_i}}.
\end{equation}

An interesting choice of weights has been designed by Floater and Hormann (FH), who introduced a family of linear barycentric rational interpolants that shows good properties in the approximation of smooth functions, in particular using equidistant nodes (see \cite{FloaterHormann}). Furthermore, in \cite{AAA}, the \emph{Adaptive Antoulas‐Anderson} (AAA) greedy algorithm for computing a barycentric rational approximant has been proposed. This scheme leads to impressively well-conditioned bases, and it has been used in computing conformal maps, or in rational minimax approximations (cf. see \cite{BDEM21} and references therein). 

Unfortunately, when the underlying function presents jump discontinuities, both the FH interpolants and the approximants produced by the AAA algorithm suffer from Gibbs effects, which can be resolved by applying the S-Gibbs algorithm adapted to this framework. Indeed, the interpolant $r_f$ admits a cardinal basis form $ r_f(x)=\sum_{j=1}^N{f_jb_j(x)}$, where
$b_j(x)= \frac{ \frac{w_j}{x-x_j}}{\sum_{i=1}^N \frac{w_i}{x-x_i}}$ is the $j$-th basis function. By composing with the map $S$, we get
$ r_f^S(x)=\sum_{i=1}^N{f_i b^S_i(x)}$, where $$b^S_j(x)= \frac{ \frac{w_j}{S(x)-S(x_j)}}{\sum_{i=1}^N \frac{w_i}{S(x)-S(x_i)}}$$ is the $j$-th mapped basis function. 
As proved in \cite{BDEM21}, by employing the S-Gibbs mapped approach we obtain a severe reduction of the Gibbs artifacts, as shown in Example \ref{eseeeempiiiio}.

\begin{example}\label{eseeeempiiiio}
Let $K_2=[-5,5]$ and
 \begin{equation*}
f_2(x)=
\begin{dcases}
    \log(-\sin(x/2))), & \text{ $-5 \le x \leq -2.5$,}\\
    \tan(x/2),& \text{ $-2.5< x \le 2$,}\\
    \arctan\left(e^{-\frac{1}{x-5.1}}\right)& \text{ $2< x \le 5$.}\\
\end{dcases}
\end{equation*}
We approximate $f_2$ both by the classical FH interpolant and by the mapped FH obtained via the S-Gibbs. In both interpolants, we used as parameter $d=8$; this parameter corresponds to the degree of the polynomial which can be reproduced exactly by using the interpolant. The results are shown in Figg. \ref{fig:figb} and \ref{fig:figbb}, where the advantages provided by the S-Gibbs scheme are evident. 
\begin{figure}[!hbt]
  \centering
  \includegraphics[width=0.45\linewidth]{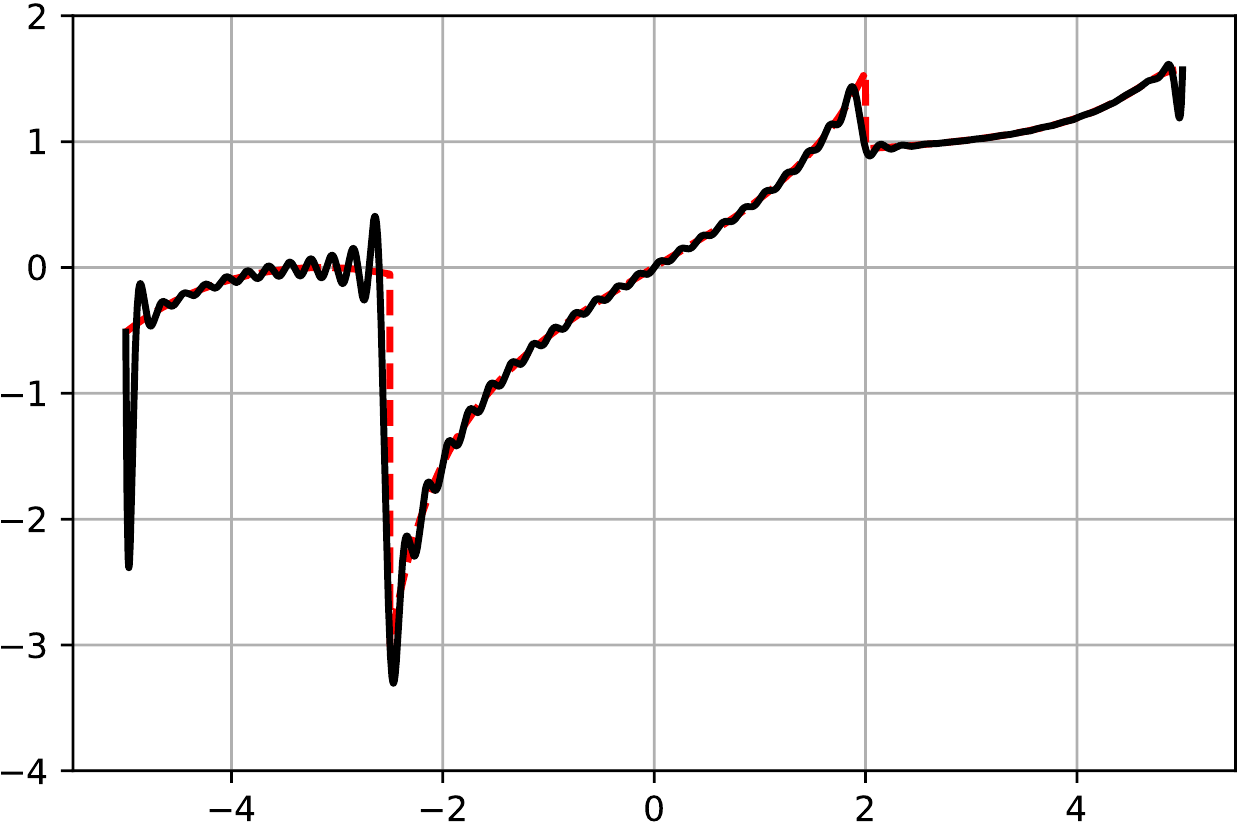} 
  \includegraphics[width=0.45\linewidth]{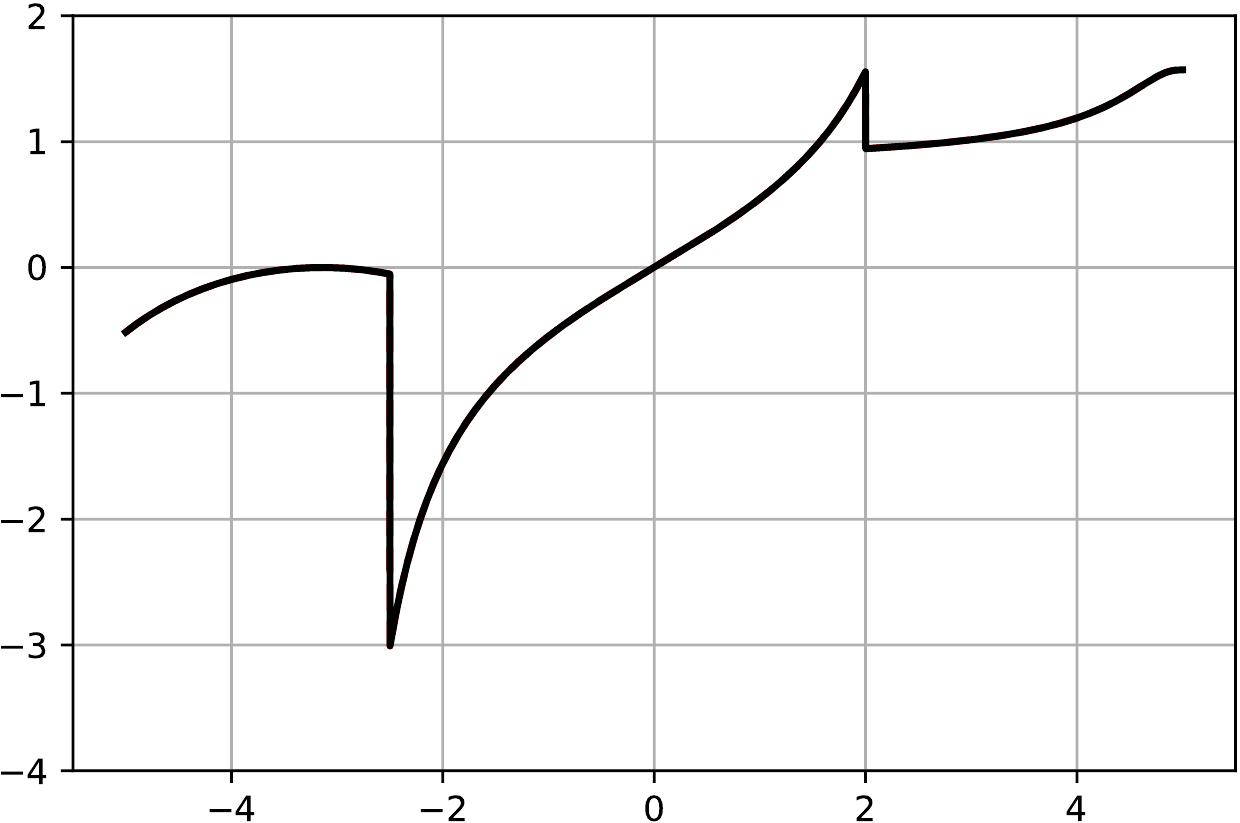}
    \caption{Approximation in $K_2$ of $f_2$. Left: the classical FH interpolant. Right: the S-Gibbs FH interpolant.}
    \label{fig:figb}
\end{figure}

\begin{figure}[!hbt]
  \centering
  \includegraphics[width=0.45\linewidth]{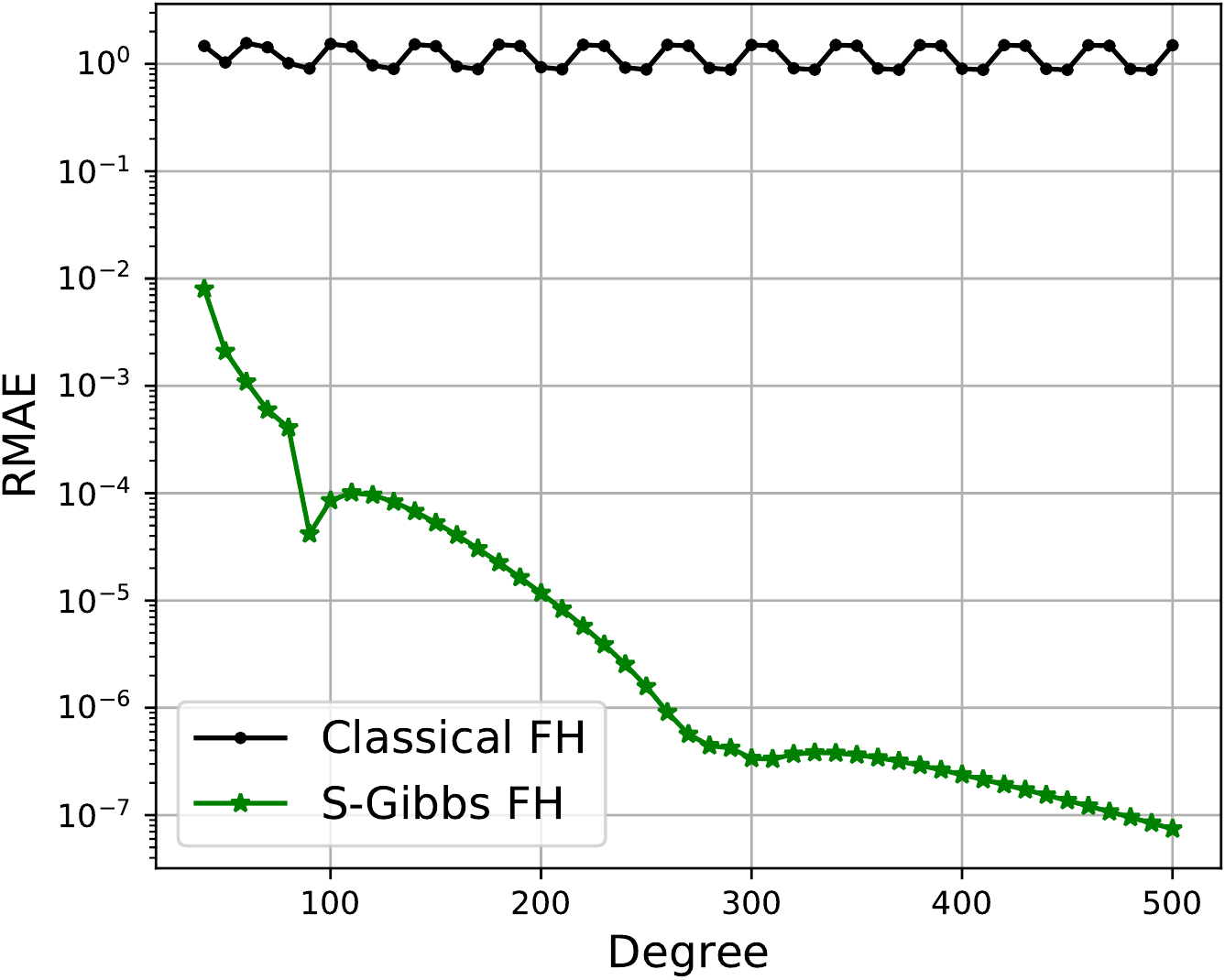}
    \caption{RMAE related to the Classical and S-Gibbs FH approaches for the function $f_2$ on $K_2$.}
    \label{fig:figbb}
\end{figure}
\end{example}

\section{Application to numerical quadrature}
Quadrature of functions at equispaced points is effective when composite rules are used. However, in the case of discontinuous functions, such rules turn out to be numerically unstable. In \cite{DeMEPP21}, first we showed that by computing analytically the quadrature weights of the \emph{mapped CL nodes}, they correspond to the quadrature weights of the trapezoidal rule. Moreover, we investigated the use of the S-Gibbs algorithm in this setting; it provides a significant improvement in the accuracy of the numerical quadrature when dealing with discontinuous functions, as shown in the following Example \ref{eeeeeesempio}. 

\begin{example}\label{eeeeeesempio}
Letting $K_3=[-2,2]$ and
 \begin{equation*}
f_3(x)=
\begin{dcases}
    \sin^2(x)-2, & \text{for $-2 \le x \leq 0.2$,}\\
\log (x^2+2)+\cos(x), & \text{for $0.2< x \le 2$,}\\
\end{dcases}. 
\end{equation*}
In Figure \ref{fig:figc}, we show the RMAE of the integral with the Newton-Cotes formulae at equispaced points, CL nodes and the quadrature formula derived via the mapped basis with the S-Gibbs scheme. Encoding the discontinuity directly into the basis via the S-Gibbs map leads to a truly effective method, which outperforms the other approaches.

\begin{figure}[!hbt]
    \centering
    \includegraphics[width=0.45\linewidth]{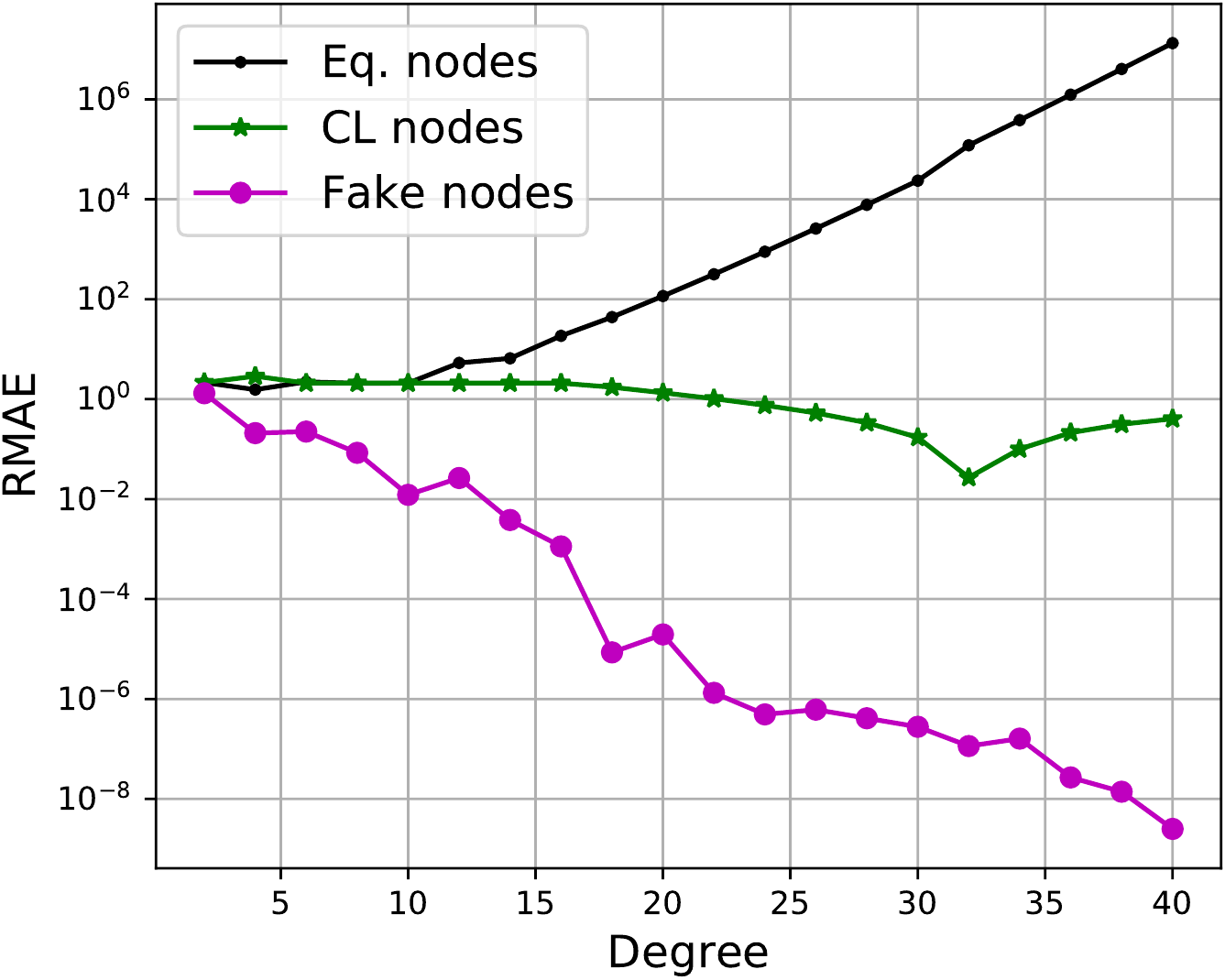}
    \caption{Approximation of the integral over $K_3$ of the functions $f_3$.}
    \label{fig:figc}
\end{figure}
\end{example}

\section{Some open problems}
\begin{itemize}
\item In the univariate setting, S-Runge and S-Gibbs have been combined in \cite{DeEM21} via the GRASPA scheme. An extension of GRASPA, at least to two dimensions, is needed.
\item GRASPA should be investigated in the framework of numerical quadrature.
\item Improved tight Lebesgue constant bounds should be investigated.
\item Recently, two dimensional mock-Chebyshev points plus regression have been investigated \cite{DellAccio22}. Is this approach a possible alternative to the mapped bases one?  
\item In multimodal medical imaging we have tested the mapped bases approach showing that indeed it is able to reduce the Gibbs artifacts \cite{PCD22}. But images are special 3-dimensional domains. A challenging problem is to extend the techniques on general $d$-dimensional domains.
\end{itemize}

\section{Conclusions}
Summarizing, the mapped bases approach proved to be effective in different approximation frameworks. Nevertheless, various related aspects ought to be explored deeper, showing that this technique still represents an active research line. Interested readers can refer to original paper \cite{4} and the ones reported in the bibliography for more insights. 

\vskip 0.1in
{\bf Acknowledgments.} This research has been accomplished within Rete ITaliana di Approssimazione (RITA), the thematic group on \lq\lq Approximation Theory and Applications\rq\rq (TAA) of the Italian Mathematical Union and partially funded by GNCS-IN$\delta$AM. We are very grateful to the former collaborators of this project: Wolfgang Erb (University of Padova), Emma Perracchione (Polytechnic of Torino) and Davide Poggiali (Padova Neuroscience Center).

\end{document}